\documentclass[psreqno,reqno,11pt]{amsart}
\usepackage{eurosym}
\usepackage{amssymb}
\usepackage{tikz}
\usepackage{amsmath}
\usepackage{tikz}
\usepackage{pgflibraryplotmarks}
\usepackage{wasysym}
\usepackage{stmaryrd}

\theoremstyle{plain}

\theoremstyle{definition}

\theoremstyle{remark}

\begin{document}
\def\sect#1{\section*{\leftline{\large\bf #1}}}
\def\th#1{\noindent{\bf #1}\bgroup\it}
\def\endth{\egroup\par}

\title[Functional Completions]{%
Functional Completions of Archimedean Vector Lattices}
\author{G. Buskes}
\address{Department of Mathematics, University of Mississippi, University,
MS 38677, USA}
\email{mmbuskes@olemiss.edu}
\author{C. Schwanke}
\address{Department of Mathematics, University of Mississippi, University,
MS 38677, USA}
\email{cmschwan@olemiss.edu}
\date{\today}

\begin{abstract}
We study completions of Archimedean vector lattices relative to any nonempty set of positively-homogeneous functions on finite-dimensional real vector spaces. Examples of such completions include square mean closed and geometric closed vector lattices, amongst others. These functional completions also lead to a universal definition of the complexification of any Archimedean vector lattice and a theory of tensor products and powers of complex vector lattices in a companion paper.
\end{abstract}

\maketitle

\vspace{.2in}
\noindent
\textbf{1. Introduction}
\vspace{.2in}

The current paper lays the foundation for a theory of polynomials on complex vector lattices. Since recent development of polynomials on real vector lattices have required a simultaneous use and development of the Fremlin tensor product (see \cite{BuBus}, \cite{BuBusPopTcacTroit}, \cite{BusvR2}), we first started to investigate the Fremlin tensor product $C(X)\bar{\otimes}C(Y)$ for compact Hausdorff spaces $X$ and $Y$ and found in Theorem 4.11 of \cite{BusSch2} that $C(X)\bar{\otimes}C(Y)+i(C(X)\bar{\otimes}C(Y))$ often cannot be equipped with a complex modulus as studied in \cite{MW}, required for it to be a complex vector lattice. Construction of a completion to add a modulus pointed itself as a natural way to proceed. Surprisingly, though the use of complex vector lattices and complex Banach lattices is more than half a century old, no description of complex vector lattices as a completion exists in the literature. In fact, a development of the theory of complex vector lattices has suffered from what appears an almost universal blanket assumption of uniform completeness in order to have a modulus available.

As we studied an appropriate completion, we found a host of similar completions, which we call functional completions, that are equally useful in applications to existing literature in vector lattice theory. Indeed, for any continuous, positively-homogeneous, real-valued function $h$ on $\mathbb{R}^{n}$ one can define a functional calculus on any uniformly complete Riesz space. The smallest vector lattice in the uniform completion of an Archimedean vector lattice $E$ on which such a calculus can be defined will be called the $h$-completion of $E$. One of these completions, associated with $h(x,y)=\sqrt{x^{2}+y^{2}}$ is needed for the completion connected with complex vector lattices mentioned above. However, the material that deals with complex vector lattices, complexification of Archimedean vector lattices, and Fremlin tensor products of complex vector lattices gradually separated itself off into a separate companion paper \cite{BusSch2}.

In the current paper, we define the functional completions universally and for their construction we need the uniform completion. However, we were unable to find in the literature the exact construction of the uniform completion of an Archimedean vector lattice that is most useful for this paper. Thus we develop the theory of uniform completions (Proposition 3.2) as well and construct functional completions relative to any nonempty set of continuous, positively-homogeneous, real-valued functions within these uniform completions (Theorem 3.18). Our uniform completion bridges (via an idea in a paper by Triki \cite{Tri}) Quinn's uniform completion in \cite{Qui} with Veksler's uniform completion in \cite{Vek}.

The theory of functional completions itself clarifies and extends previous results in the literature for very specific positively-homogeneous functions like the square mean and the geometric mean (see, e.g., \cite{AzBoBus}, Corollary 4.7, \cite{Az}, Proposition 2.20, and Corollaries 3.9 and 3.10 in this paper) and forms a foundation for our treatment of complex vector lattices in \cite{BusSch2}.

Finally, we connect the use of differential calculus as first seen in Theorem 4.2 of \cite{BeuHuidP} by Beukers, Huijsmans, and de Pagter to $h$-completions for convex or concave $h$. By doing so in Theorem 3.8, we sharpen a special case of Kusraev's Theorem 5.5 in \cite{Kus} (while keeping the structure of its proof largely intact) in three ways. We weaken the assumption of uniform completeness, verify that the proof of Theorem 5.5 in \cite{Kus} in our special case does not (contrary to Kusraev's proof) require more than the Countable Axiom of Choice, and provide more concrete formulas that directly link to Lemmas 4.2 and 4.3 in \cite{AzBoBus}. As indicated in the previous sentence, we keep this paper in the framework of \cite{BusdPvR} and \cite{BusvR3} in not using more than a modicum of the Axiom of Choice.

\vspace{.5in}
\noindent
\textbf{2. Preliminaries}
\vspace{.2in}

For all unexplained terminology about vector lattices we refer the reader to
the standard texts \cite{AB}, \cite{LuxZan1}, and \cite{Zan2}. Throughout, $\mathbb{R}$ stands for the set of real numbers while $\mathbb{N}$ stands for the set of (nonzero) positive integers. For $s\in\mathbb{N}$ and sets $A_{1},...,A_{s}$, we write $\times_{k=1}^{s}A_{k}$ for the Cartesian product $A_{1}\times\dots\times A_{s}$. In case that $A_{k}=A$ for every $k\in\lbrace 1,...,s\rbrace$ we write $\times_{k=1}^{s}A_{k}=A^{s}$. For an Archimedean vector lattice $E$, we denote the \textit{positive cone} $\lbrace f\in E:f\geq 0\rbrace$ of $E$ by $E^{+}$ and we denote the Dedekind completion of $E$ by $E^{\delta}$.

Let $E_{1},...,E_{s},F$ be Archimedean vector lattices. We say that a map $T:\times_{k=1}^{s}E_{k}\rightarrow F$ that is linear in each variable separately is \textit{$s$-linear}. An $s$-linear map $T:\times_{k=1}^{s}E_{k}\rightarrow F$ is called a \textit{vector lattice $s$-morphism} if the map $f_{k}\mapsto T(f_{1},...,f_{k},...,f_{s})\ (f_{k}\in E_{k})$ is a vector lattice homomorphism for each $k\in\lbrace 1,...,s\rbrace$ and for all $f_{j}\in E_{j}^{+}\ (j\neq k)$.

We will freely use functional calculus for Archimedean vector
lattices (see \cite{BusdPvR}). We write
$\mathcal{H}(\mathbb{R}^{m})$ for the space of all continuous, real-valued functions $h$ on $\mathbb{R}^{m}$ that are \textit{positively-homogeneous}, i.e. $h(\lambda x)=\lambda h(x)$ for every $\lambda
\in \mathbb{R}^{+}$ and all $x\in \mathbb{R}^{m}$. The space of all nonzero real-valued vector lattice homomorphisms on an Archimedean vector lattice $E$ is denoted by $H(E)$. For an Archimedean vector lattice $E$ and a nonempty subset $A$ of $E$, we denote by $[A]$ the vector space generated by $A$ in $E$ and we denote by $\langle A\rangle$ the vector sublattice generated by $A$ in $E$. For an Archimedean vector lattice $E$, with $a_{1},...,a_{m},b\in E$ and $h\in\mathcal{H}(\mathbb{R}^{m})$, we write $h(a_{1},...,a_{m})=b$ when $h(\omega
(a_{1}),...,\omega (a_{m}))=\omega (b)$ for every $\omega\in H(\left\langle a_{1},...,a_{m},b\right\rangle)$.

We remind the reader of the definition of the uniform completion of an Archimedean vector lattice.

\smallskip
\noindent
\th{Definition 2.1.} Given an
Archimedean vector lattice $E$, a
sequence $(f_{n})$ in $E$ is said to converge relatively
uniformly to $f$ in $E$ if there exists $0<p\in E$
such that for every $\epsilon>0$ there exists $N\in\mathbb{N}$ for which $|f_{n}-f|<\epsilon p$ for every $n\geq N$. In this case, we write $f_{n}\overset{ru}{\rightarrow }f$. We call a sequence $(f_{n})$ in $E$ a relatively uniformly Cauchy sequence if there exists $0<p\in E$ such that for every $\epsilon >0$ there exists $N\in\mathbb{N}$ for which $|f_{m}-f_{n}|<\epsilon p$ for every $m,n\geq N$. If every relatively uniformly Cauchy sequence in $E$ converges relatively uniformly in $E$, we say that $E$ is uniformly complete.
\endth

\smallskip
Note that in the above, there exists at most one $f$ such that $f_{n}\overset{ru}{\rightarrow }f$ (see Theorem 16.2(i) in \cite{LuxZan1}). There exist various ways of introducing uniform completions
of Archimedean vector lattices in the literature, (see \cite{vH}, \cite{Qui}, and \cite{Vek}). For our purposes, we choose the definition by van Haandel in \cite{vH}.

\smallskip
\noindent
\th{Definition 2.2.}\ \textnormal{(\cite{vH}, Definition 8.6)}\ Given an
Archimedean vector lattice $E$, we call
a pair $(E^{u},\phi)$ a uniform completion of $E$ if
the following hold.
\begin{itemize}
\item[(1)] $E^{u}$ is a uniformly complete
Archimedean vector lattice.
\item[(2)] $\phi:E\rightarrow E^{u}$ is an injective vector lattice homomorphism.
\item[(3)] For every uniformly complete Archimedean vector lattice $F$ and for every vector lattice homomorphism $T:E\rightarrow F$ there exists a unique vector lattice homomorphism $T^{u}:E^{u}\rightarrow F$ such that $T^{u}\circ\phi=T$.
\end{itemize}
\endth

\smallskip
We will also use the following definition, which was introduced (with slightly
different notation) on page 85 of \cite{LuxZan1}. For an Archimedean vector lattice $E$ and for $A\subseteq E$, we define the \textit{pseudo uniform closure} $\bar{A}$ of $A$ to be the set of all $f\in E$ for which there exists a sequence $(f_{n})$ in $A$ such
that $f_{n}\overset{ru}{\rightarrow }f$. We call $A$ \textit{relatively uniformly closed} if $\bar{A}=A$. The relatively uniformly closed sets are the closed sets in the \textit{relatively uniform topology}, defined in \cite{LuxZan1}.

Finally, we iterate the pseudo uniform closure of a nonempty subset $L$ of $E$ as follows via transfinite induction.
\begin{equation*}
\begin{array}{l}
A_{1}:=A\text{,} \\ 
A_{\alpha}:=\overline{A_{\alpha-1}}\ \text{when}\ \alpha>1\ \text{is not a limit ordinal, and}\\ 
A_{\alpha}:=\bigcup_{\beta<\alpha}A_{\beta}\ \text{when}\ \alpha\ \text{is a limit ordinal.}
\end{array}
\end{equation*}

\vspace{.5in}
\noindent
\textbf{3. Completions}
\vspace{.2in}

In this section we select certain intermediate vector lattices (see \cite{Qui}) as completions of Archimedean vector lattices
via nonempty subsets of $\bigcup_{m\in \mathbb{N}}\mathcal{H}(\mathbb{R}^{m})$. In \cite{BusSch2}, we
employ one of these completions to complexify Archimedean vector lattices. Since \cite{vH} is somewhat inaccessible and the proof of the existence of the uniform completion in \cite{vH} skips the use of the iterated pseudo-closures, we provide a different proof. We start by extending positive linear maps on vector sublattices of an Archimedean vector lattice to their pseudo-closures as follows.

\smallskip
\noindent
\th{Proposition 3.1.} Let $L$ be a vector
sublattice of an Archimedean vector lattice $E$. Then the
following hold.
\begin{itemize}
\item[(1)] $L_{\alpha}$ is a vector sublattice of $E$ for every
ordinal $\alpha$.
\item[(2)] $L_{\omega_{1}}$ is relatively uniformly closed in $E$.
\item[(3)] If $E$ is relatively uniformly complete then so is $L_{\omega_{1}}$.
\item[(4)] $L$ is dense in $L_{\omega_{1}}$ in the relatively uniform topology.
\item[(5)] For every uniformly complete Archimedean vector lattice $F$, every ordinal $1\leq\alpha\leq\omega_{1}$, and every positive linear map $T:L\rightarrow F$ there exists a unique positive linear map $T_{\alpha}:L_{\alpha}\rightarrow F$ such that $T_{\alpha}|_{L}=T$. Moreover, if $T$ is a vector lattice homomorphism then so is $T_{\alpha}$.
\end{itemize}
\endth

\begin{proof} Statement (1) follows from transfinite induction and uses elementary calculus of relatively uniformly convergent sequences, (see Theorem 16.2 of \cite{LuxZan1}). Part (2) is an immediate consequence of the fact that every sequence in $L_{\omega_{1}}$ resides in an $L_{\alpha}$ for some $\alpha<\omega_{1}$, and (3) follows directly from (2), whereas (4) follows directly from the definition of the relatively uniform topology. To prove (5), let $T:L\rightarrow F$ be a positive linear map and define $T_{1}:=T$. Next, let $1<\alpha\leq\omega_{1}$ be an ordinal and assume that $T$ can be uniquely extended to a positive linear map $T_{\beta}:L_{\beta}\rightarrow F$ for every ordinal $1\leq\beta<\alpha$. Let $f\in L_{\alpha}$. Suppose that $\alpha$ is not a limit ordinal. There exists a relatively uniformly Cauchy sequence $(f_{n})$ in $L_{\alpha-1}$ such that $f_{n}\overset{ru}{\rightarrow}f$. Since $T_{\alpha-1}$ is positive,
\begin{equation*}
|T_{\alpha-1}(g)|\leq T_{\alpha-1}(|g|)\ \text{for all}\ g\in L_{\alpha-1}.\tag{*}
\end{equation*}
\noindent
Therefore $(T_{\alpha-1}(f_{n}))$ is a relatively uniformly Cauchy sequence in the uniformly complete $F$. Hence there exists (a unique) $h\in F$ such that $T_{\alpha-1}(f_{n})\overset{ru}{\rightarrow}h$. Define $T_{\alpha}:L_{\alpha}\rightarrow F$ by $T_{\alpha}(f)=h$. It follows from (*) that $T_{\alpha}$ is well-defined. If $\alpha$ is a limit ordinal then define $T_{\alpha}(f)=T_{\beta}(f)\ (f\in L_{\beta}\ \text{and}\ \beta<\alpha)$. By the induction hypothesis, $T_{\alpha}$ is well-defined.

It is readily checked by using elementary calculus of relatively uniformly convergent sequences that $T_{\alpha}$ is a positive linear map for every ordinal $1\leq\alpha\leq\omega_{1}$, and that $T_{\alpha}$ is a 
vector lattice homomorphism if $T$ is a vector lattice homomorphism. That $T_{\alpha}$ is indeed the unique positive linear extension of $T$ to $L_{\alpha}$ follows from uniform density and transfinite induction.
\end{proof}

\smallskip
It is evident that a uniform completion, if it exists, is unique. We use the previous proposition to prove that every Archimedean vector lattice has a uniform completion. The reader should compare Proposition 3.2 with Theorem 3.3 of \cite{Tri}, where Triki deals with Quinn's definition of uniform completion (see \cite{Qui}). A small adaptation of Theorem 3.3 of \cite{Tri} to vector lattice homomorphisms rather than positive linear maps shows, in effect, that Quinn's definition of uniform completion is equivalent to van Haandel's definition above. In addition, we now generalize Theorem 3.3 of \cite{Tri} to multilinear maps.

\smallskip
\noindent
\th{Proposition 3.2.} \textnormal{(1)} If $E$ is an Archimedean vector lattice then there
exists a uniform completion $(E^{u},\phi)$ of $E$.
\begin{itemize}
\item[(2)] If $E_{1},...,E_{s},F$ are Archimedean vector lattices with $F$ uniformly complete and $T:\times_{k=1}^{s}E_{k}\rightarrow F$ is positive and $s$-linear then there exist injective vector lattice homomorphisms $\phi_{k}:E_{k}\rightarrow E_{k}^{u}\ (k\in\lbrace 1,...s\rbrace)$ as well as a unique positive $s$-linear map $T^{u}:\times_{k=1}^{s}E_{k}^{u}\rightarrow F$ such that $T^{u}(\phi_{1}(f_{1}),...,\phi_{s}(f_{s}))=T(f_{1},...,f_{s})$ for every $(f_{1},...,f_{s})\in \times _{k=1}^{s}E_{k}$. Moreover, if $T$ is a vector lattice $s$-morphism then so is $T^{u}$.
\end{itemize}
\endth

\begin{proof}
(1) Assume that $E$ and $F$ are Archimedean vector lattices. The natural embedding $\phi$ of $E$ into $E^{\delta}$ yields an injective vector lattice homomorphism. Define
\begin{equation*}
E^{u}:=\phi(E)_{\omega_{1}}.
\end{equation*}
\noindent
Since $\phi(E)\subseteq E^{u}$, we may consider $\phi$ as a map from $E$ to $E^{u}$. Let $T:E\rightarrow F$ be a positive linear map. Then the map $\tilde{T}:\phi(E)\rightarrow F$ defined by $\tilde{T}(\phi(f))=T(f)$ is also a positive linear map, and if $T$ is a vector lattice homomorphism then so $\tilde{T}$. By (5) of Proposition 3.1, there exists a unique positive linear extension $\tilde{T}_{\omega_{1}}:E^{u}\rightarrow F$ of $\tilde{T}$, and if $\tilde{T}$ is a vector lattice homomorphism then so is $\tilde{T}_{\omega_{1}}$. Also, $\tilde{T}_{\omega_{1}}\circ\phi=T$.

\smallskip
(2) Let $E_{1},...,E_{s},F$ be Archimedean vector lattices. For each $k\in\lbrace 1,...,s\rbrace$, let $\phi_{k}$ be the natural embedding of $E$ into $E^{\delta}$, considered as a map from $E_{k}$ to $\phi_{k}(E)$. Define
\begin{equation*}
E_{k}^{u}:=\phi_{k}(E_{k})_{\omega_{1}}
\end{equation*}
\noindent
for each $k\in\lbrace 1,...,s\rbrace$. Suppose $T:\times_{k=1}^{s}E_{k}\rightarrow F$ be a positive $s$-linear map and consider 
$T$ to be a map from $\times_{k=1}^{s}\phi_{k}(E_{k})$ to $F$ by identifying $\phi_{k}(f_{k})$ with $f_{k}$ for all 
$f_{k}\in E_{k}\ (k\in\lbrace 1,...,s\rbrace)$. For every $g_{k}\in E_{k}^{+}\ (k\in\lbrace 2,...,s\rbrace)$ we define 
\begin{align*}
T_{g_{2},...,g_{s}}(x)=T(x,g_{2},...,g_{s})\ (x\in E_{1}).
\end{align*}
\noindent
By (5) in Proposition 3.1, there exists a unique positive linear map $T_{g_{2},...,g_{s}}^{u}:E_{1}^{u}\rightarrow F$ such that $T_{g_{2},...,g_{s}}^{u}(x)=T_{g_{2},...,g_{s}}(x)\ (x\in E_{1})$. Moreover, if $T_{g_{2},...,g_{s}}$ is a vector lattice homomorphism then so is $T_{g_{2},...,g_{s}}^{u}$. Next, define
\begin{align*}
T_{+}(g_{1},...,g_{s})=T^{u}_{g_{2},...,g_{s}}(g_{1})\ (g_{1}\in E_{1}^{u}\ \text{and}\ g_{k}\in E_{k}^{+}\ (k\in\lbrace 2,...,s\rbrace))
\end{align*}
Let $j\in\lbrace 2,...,s\rbrace$ and let $g_{j},g'_{j}\in E_{j}^{+}$. Since $T^{u}_{g_{2},...g_{j}+g'_{j},...,g_{s}}$ and $T^{u}_{g_{2},...g_{j},...,g_{s}}+T^{u}_{g_{2},...g'_{j},...,g_{s}}$ are both positive linear extensions of $T_{g_{2},...,g_{j}+g'_{j},...,g_{s}}$ from $E_{1}^{u}$ to $F$, it follows from the uniqueness of such extensions that $T^{u}_{g_{2},...g_{j}+g'_{j},...,g_{s}}=T^{u}_{g_{2},...g_{j},...,g_{s}}+T^{u}_{g_{2},...g'_{j},...,g_{s}}$. Therefore, $T_{+}$ is additive in each variable separately. By routine reasoning, $T_{+}$ extends to a positive $s$-linear map from $E_{1}^{u}\times E_{2}\times\dots\times E_{s}$ to $F$ which is a vector lattice $s$-morphism in case $T$ is a vector lattice $s$-morphism. By repeating this argument for the remaining $s-1$ variables, we obtain the desired result.
\end{proof}

\smallskip
Next, we introduce completions of Archimedean vector lattices that are induced by nonempty subsets of $\bigcup_{m\in\mathbb{N}}\mathcal{H}(\mathbb{R}^{m})$.

\smallskip
\noindent
\th{Definition 3.3.} For an Archimedean vector lattice $E$ and for $h\in\mathcal{H}(\mathbb{R}^{m})$,
we say that $E$ is $h$-complete if for every $a_{1},...,a_{m}\in E$ there exists $b\in E$ such that $h(a_{1},...,a_{m})=b$. For a subset $\mathcal{D}$ of $\bigcup_{m\in\mathbb{N}}\mathcal{H}(\mathbb{R}^{m})$, we say that $E$ is $\mathcal{D}$-complete if $E$ is $h$-complete for every $h\in\mathcal{D}$.
\endth

\smallskip
The Stolarsky and Gini means (see \cite{Sto}, respectively \cite{NeuPal}) are well-studied examples of elements of $\mathcal{H}(\mathbb{R}^{2})$. Though they are typically defined on $(0,\infty)^{2}$, but can be extended continuously to all of $\mathbb{R}^{2}$ as follows.

\smallskip
\noindent
\th{Example 3.4.} For real numbers $r\neq s$ and $s\neq 0$, define
\begin{align*}
\mu_{r,s}(x,y)=\left\{
     \begin{array}{lr}
       \bigl(\frac{r(|x|^{s}-|y|^{s})}{s(|x|^{r}-|y|^{r})}\bigr)^{\frac{1}{s-r}} & if\ x\neq y\\
       |x| & if\ x=y
     \end{array}
   \right.
\end{align*}
\noindent
for $x,y\in\mathbb{R}$. We call $\mu_{r,s}$ the $(r,s)$-Stolarsky mean. Particularly, $\mu_{2,4}(x,y)=\sqrt{\frac{|x|^{2}+|y|^{2}}{2}}$ for $x,y\in\mathbb{R}$ and $\mu_{1,-1}(x,y)=\sqrt{|xy|}$ for $x,y\in\mathbb{R}$. We call $\mu_{2,4}$ the square mean and $\mu_{1,-1}$ the geometric mean.
\endth

\smallskip
\noindent
\th{Example 3.5.}  For real numbers $r\neq s$, define
\begin{align*}
\nu_{r,s}(x,y)=\left\{
     \begin{array}{lr}
       \bigl(\frac{|x|^{s}+|y|^{s}}{|x|^{r}+|y|^{r}}\bigr)^{\frac{1}{s-r}} & if\ (x,y)\neq (0,0)\\
       0 & if\ (x,y)=(0,0)
     \end{array}
   \right.
\end{align*}
\noindent
for $x,y\in\mathbb{R}$. We call $\nu_{r,s}$ the $(r,s)$-Gini mean.
\endth

\smallskip
In order to complexify Archimedean vector lattices via functional completions, we connect functional calculus to the modulus formula
\begin{align*}
|f+ig|=\sup\lbrace f\cos\theta+g\sin\theta:\theta\in[0,2\pi]\rbrace
\end{align*}
\noindent
for complex vector lattices, (see, e.g, Section 91 of \cite{Zan2}). Following the idea to use tangents by Beukers, Huijsmans, and de Pagter in Theorem 4.2 of \cite{BeuHuidP}, we identify elements of an Archimedean vector lattice $E$ of the form $h(a_{1},...,a_{m})\ (a_{1},...,a_{m}\in E^{+})$ for convex or concave $h\in\mathcal{H}(\mathbb{R}^{m})$ with elements of $E$ that are defined via differential calculus. To this end, we need some notations.

\smallskip
\noindent
\textbf{Notations 3.6.} Let $E$ be an Archimedean vector lattice. The Euclidean norm on $\mathbb{R}^{m}$ is denoted by $||\ ||$. For $h\in\mathcal{H}(\mathbb{R}^{m})$ we set
\begin{align*}
\Delta_{h}=\lbrace\text{\textit{\textbf{c}}}\in(\mathbb{R}^{+})^{m}:h\ \text{is differentiable at}\ \text{\textit{\textbf{c}}}\ \text{and}\ ||\text{\textit{\textbf{c}}}||=1\rbrace.
\end{align*}
\noindent
For $h\in\mathcal{H}(\mathbb{R}^{m}), \text{\textit{\textbf{c}}}\in\Delta_{h}$, and ${\text{\textit{\textbf{a}}}}:=(a_{1},...,a_{m})\in E^{m}$ we define $\nabla h(\text{\textit{\textbf{c}}})\cdot{\text{\textit{\textbf{a}}}}:=\sum\limits_{k=1}^{m}\frac{\partial h(\text{\textit{\textbf{c}}})}{\partial x_{k}}a_{k}$. For $a_{1},...,a_{m}\in E^{+}\ (m\geq 2)$, and for $\theta=(\theta_{1},...,\theta_{m-1})\in[0,\pi]^{m-1}$ we define
\begin{center}
$s_{\theta}(a_{1},...,a_{m}):=$
\end{center}
\begin{align*}
\cos\theta_{1}a_{1}+\sum\limits_{k=2}^{m-2}\Bigl(\prod\limits_{j=1}^{k-1}\sin\theta_{j}\Bigr)\cos\theta_{k}a_{k}+\Bigl(\prod\limits_{j=1}^{m-2}\sin\theta_{j}\Bigr)\cos\theta_{m-1}a_{m-1}+\Bigl(\prod\limits_{j=1}^{m-1}\sin\theta_{j}\Bigr)a_{m},
\end{align*}
\noindent
where $\sum\limits_{k=2}^{m-2}\Bigl(\prod\limits_{j=1}^{k-1}\sin\theta_{j}\Bigr)\cos\theta_{k}$ is taken to equal zero for $m\in\lbrace 2,3\rbrace$, and $\Bigl(\prod\limits_{j=1}^{m-2}\sin\theta_{j}\Bigr)\cos\theta_{m-1}$ is taken to equal zero for $m=2$. For short, we denote $s_{\theta}(e_{1},...,e_{m})$ by $s_{\theta}$, where $e_{k}$ is the $k$th element of the standard orthonormal basis of $\mathbb{R}^{m}$. Finally, for $n\in\mathbb{N}$ we also set
\begin{align*}
P_{n}=\lbrace(\frac{l_{1}\pi}{2^{n+1}},...,\frac{l_{m-1}\pi}{2^{n+1}}):l_{1},...,l_{m-1}\in\lbrace 1,...,2^{n}\rbrace\rbrace.
\end{align*}

In Theorem 3.8, we give, as described in the introduction, a sharpening of a special case of Theorem 5.5 in \cite{Kus}. We need the following lemma. Parts (1) and (2) are known, and part (3) may also be known, but we were unable to find a reference.

\smallskip
\noindent
\th{Lemma 3.7.} Let $h\in\mathcal{H}(\mathbb{R}^{m})$ be convex \textnormal{(}respectively, concave\textnormal{)} on $(\mathbb{R}^{+})^{m}$. Then
\begin{itemize}
\item[(1)] $h$ is differentiable almost everywhere with respect to Lebesgue measure on $(0,\infty)^{m}$.
\item[(2)] If $h$ is differentiable at $\textit{\textbf{x}}\in\mathbb{R}^{m}$ then $\nabla h(\lambda\textit{\textbf{x}})=\nabla h(\textit{\textbf{x}})$ for every $0<\lambda\in\mathbb{R}$.
\item[(3)] $h(\textit{\textbf{x}})=\sup\lbrace\nabla h(\text{\textit{\textbf{c}}})\cdot{\text{\textit{\textbf{x}}}}:\text{\textit{\textbf{c}}}\in\Delta_{h}\rbrace$ for every ${\text{\textit{\textbf{x}}}}\in(\mathbb{R}^{+})^{m}$, \textnormal{(}respectively, $h(\textit{\textbf{x}})=\inf\lbrace\nabla h(\text{\textit{\textbf{c}}})\cdot{\text{\textit{\textbf{x}}}}:\text{\textit{\textbf{c}}}\in\Delta_{h}\rbrace$ for every ${\text{\textit{\textbf{x}}}}\in(\mathbb{R}^{+})^{m}$\textnormal{)}.
\end{itemize}
\endth

\begin{proof} (1) By Exercise 1.17 in \cite{Phel}, which follows from Radamacher's Theorem (also see Exercise 1.18 in \cite{Phel}), $h$ is differentiable on $(0,\infty)^{m}$ outside a set of Lebesgue measure zero.

(2) Follows directly from the positive homogeneity of $h$.

(3) Suppose that $h$  is convex on $(\mathbb{R}^{+})^{m}$. It follows from Euler's Homogeneous Function Theorem (for instance, Exercise 2-34 in \cite{Spiv}) that $\nabla h(\text{\textit{\textbf{c}}})\cdot\text{\textit{\textbf{c}}}=h(\text{\textit{\textbf{c}}})$ whenever $\text{\textit{\textbf{c}}}\in\mathbb{R}^{m}$ and $h$ is differentiable at $\text{\textit{\textbf{c}}}$. From this observation as well as the convexity of $h$, we obtain
\begin{align*}
h(\textit{\textbf{x}})=\sup\lbrace\nabla h(\text{\textit{\textbf{c}}})\cdot{\text{\textit{\textbf{x}}}}:\text{\textit{\textbf{c}}}\in\Delta_{h}\rbrace
\end{align*}
\noindent
for every ${\text{\textit{\textbf{x}}}}\in(\mathbb{R}^{+})^{m}$ where $h$ is differentiable, as well as $\nabla h(\text{\textit{\textbf{c}}})\cdot\text{\textit{\textbf{x}}}\leq h(\text{\textit{\textbf{x}}})$ for every $\text{\textit{\textbf{c}}}\in\Delta_{h}$ and for every $\text{\textit{\textbf{x}}}\in(\mathbb{R}^{+})^{m}$. Suppose that $h$ is not differentiable at $\text{\textit{\textbf{b}}}\in(\mathbb{R}^{+})^{m}$, and let $\epsilon>0$. We just need to show that there exists $\text{\textit{\textbf{c}}}\in\Delta_{h}$ such that $h(\text{\textit{\textbf{b}}})-\nabla h(\text{\textit{\textbf{c}}})\cdot\text{\textit{\textbf{b}}}<\epsilon$. To this end, we note that since $h$ is continuous, there exists $\delta_{1}>0$ such that $||h(\text{\textit{\textbf{c}}})-h(\text{\textit{\textbf{b}}})||<\epsilon/2$ whenever $||\text{\textit{\textbf{c}}}-\text{\textit{\textbf{b}}}||<\delta_{1}$. Since $h$ is convex and continuous, it is locally Lipschitz (see \cite{Phel}, Proposition 1.6). Let $B_{\delta_{2}}(\text{\textit{\textbf{b}}})=\lbrace\text{\textit{\textbf{x}}}\in(\mathbb{R}^{+})^{m}:||\text{\textit{\textbf{x}}}-\text{\textit{\textbf{b}}}||<\delta_{2}\rbrace$ be a neighborhood of $\text{\textit{\textbf{b}}}$ where $h$ is Lipschitz, say with Lipschitz constant $M$. Then each partial derivative satisfies $|\frac{\partial h(\text{\textit{\textbf{c}}})}{\partial x_{k}}|\leq M$ for every $\text{\textit{\textbf{c}}}\in B_{\delta_{2}}(\text{\textit{\textbf{b}}})$ where $h$ is differentiable. Furthermore, there exists $\delta_{3}>0$ such that if $||\text{\textit{\textbf{c}}}-\text{\textit{\textbf{b}}}||<\delta_{3}$ then $\sum_{k=1}^{m}|c_{k}-c_{x_{k}}|<\frac{\epsilon}{2mM}$. In this case we have $||\nabla h(\text{\textit{\textbf{c}}})\cdot \text{\textit{\textbf{c}}}-\nabla h(\text{\textit{\textbf{c}}})\cdot \text{\textit{\textbf{b}}}||=||\nabla h(\text{\textit{\textbf{c}}})\cdot(\text{\textit{\textbf{c}}}-\text{\textit{\textbf{b}}})||\leq mM\sum_{k=1}^{m}|c_{k}-c_{x_{k}}|<\epsilon/2$. Since $h$ is differentiable almost everywhere on $(0,\infty)^{m}$, we may choose $\text{\textit{\textbf{c}}}\in\Delta_{h}$ and $0<\lambda\in\mathbb{R}$ such that $||\lambda\text{\textit{\textbf{c}}}-\text{\textit{\textbf{b}}}||<\delta_{1}\wedge\delta_{2}\wedge\delta_{3}$. Then using part (2) and the identity $\nabla h(\text{\textit{\textbf{c}}})\cdot\text{\textit{\textbf{c}}}=h(\text{\textit{\textbf{c}}})$ for $\text{\textit{\textbf{c}}}\in\mathbb{R}^{m}$ where $h$ is differentiable,
\begin{align*}
||h(\text{\textit{\textbf{b}}})-\nabla h(\text{\textit{\textbf{c}}})\cdot \text{\textit{\textbf{b}}}||&=||h(\text{\textit{\textbf{b}}})-h(\lambda\text{\textit{\textbf{c}}})+\nabla h(\text{\textit{\textbf{c}}})\cdot \lambda\text{\textit{\textbf{c}}}-\nabla h(\text{\textit{\textbf{c}}})\cdot\text{\textit{\textbf{b}}}||\\
&\leq||h(\text{\textit{\textbf{b}}})-h(\lambda\text{\textit{\textbf{c}}})||+||\nabla h(\text{\textit{\textbf{c}}})\cdot \lambda\text{\textit{\textbf{c}}}-\nabla h(\text{\textit{\textbf{c}}})\cdot \text{\textit{\textbf{b}}}||\\
&<\epsilon/2+\epsilon/2=\epsilon.
\end{align*}
\noindent
The case where $h$ is concave on $(\mathbb{R}^{+})^{m}$ is handled in a similar manner.
\end{proof}

We now turn to the result promised before Lemma 3.7.

\smallskip
\noindent
\th{Theorem 3.8.} Let $E$ be an Archimedean vector lattice, let $h\in\mathcal{H}(\mathbb{R}^{m})$, and let $a_{1},...,a_{m}\in E^{+}$.
\begin{itemize}
\item[(1)] If $h$ is convex \textnormal{(}respectively, concave\textnormal{)} on $(\mathbb{R}^{+})^{m}$ then $h(a_{1},...,a_{m})=b$ for some $b\in E$ if and only if $b=E\text{-}\sup\lbrace\nabla h(\text{\textit{\textbf{c}}})\cdot{\text{\textit{\textbf{a}}}}:\text{\textit{\textbf{c}}}\in\Delta_{h}\rbrace$ \textnormal{(}respectively, $b=E\text{-}\inf\lbrace\nabla h(\text{\textit{\textbf{c}}})\cdot{\text{\textit{\textbf{a}}}}:\text{\textit{\textbf{c}}}\in\Delta_{h}\rbrace$\textnormal{)}.
\item[(2)] If $h$ is convex (respectively, concave), and if $m\geq 2$, and if all the partial derivatives of $h$ are uniformly continuous on $\lbrace s_{\theta}:\theta\in\bigcup_{n\in\mathbb{N}}P_{n}\rbrace$ then the sequence $\bigl(\sup\lbrace\nabla h(s_{\theta})\cdot{\text{\textit{\textbf{a}}}}:\theta\in P_{n}\rbrace\bigr)$ \textnormal{(}respectively, the sequence $\bigl(\inf\lbrace\nabla h(s_{\theta})\cdot{\text{\textit{\textbf{a}}}}:\theta\in P_{n}\rbrace\bigr)$\textnormal{)} converges relatively uniformly to $h(a_{1},...,a_{m})$.
\end{itemize}
\endth

\begin{proof} (1) Let ${\text{\textit{\textbf{a}}}}=(a_{1},...,a_{m})\in E^{m}$ and write $A=\lbrace\nabla h(\text{\textit{\textbf{c}}})\cdot{\text{\textit{\textbf{a}}}}:\text{\textit{\textbf{c}}}\in\Delta_{h}\rbrace$. Suppose that $E$-$\sup A$ exists in $E$. To prove that $E$-$\sup A=h(a_{1},...,a_{m})$ we need to prove that $h(\omega(a_{1}),...,\omega(a_{m}))=\omega(E\text{-}\sup A)$ for every $\omega$ in a point separating subset of $H(L)$, where $L:=\langle a_{1},...,a_{m}, E\text{-}\sup A\rangle$ (see Lemma 3.3 of \cite{BusdPvR}). To this end, we note that there exists a metrizable space $Y$, an order-dense vector sublattice $F$ of $C(Y)$, and a vector lattice isomorphism $\phi:L\rightarrow F$ (see \cite{BusvR3}, (ii) on page 526 and Theorem 2.4(i)). Consequently, we can take the point separating set of the previous sentence to be the set of all compositions of point evaluations $\hat{y}\ (y\in Y)$ with $\phi$. Moreover, $F$-$\sup\phi(A)=C(Y)$-$\sup\phi(A)$ (see \cite{dJvR}, Lemma 13.21(i)). Since $\phi$ is an isomorphism and since $E$-$\sup A\in L$ we have
\begin{align*}
\phi(E\text{-}\sup A)=\phi(L\text{-}\sup A)=F\text{-}\sup\phi(A)=C(Y)\text{-}\sup\phi(A).
\end{align*}
Define $b(y)=h(\phi(a_{1})(y),...,\phi(a_{m})(y))\ (y\in Y)$. Then $b\in C(Y)$. Since $a_{1},...,a_{m}\in E^{+}$ we have $\phi(a_{k})\in C(Y)^{+}\ (k\in\lbrace 1,...,m\rbrace)$. From Lemma 3.7(3) above,
\begin{align*}
b(y)=\sup\Bigl\lbrace\sum\limits_{k=1}^{m}\frac{\partial h(\text{\textit{\textbf{c}}})}{\partial x_{k}}\phi(a_{k})(y):\text{\textit{\textbf{c}}}\in\Delta_{h}\Bigr\rbrace
\end{align*}
for every $y\in Y$. Therefore $b=C(Y)\text{-}\sup\phi(A)$, and thus $b=\phi(E\text{-}\sup A)$. Moreover,
\begin{align*}
\hat{y}(\phi(E\text{-}\sup A))=b(y)=h(\phi(a_{1})(y),...,\phi(a_{m})(y))=h(\hat{y}(\phi(a_{1})),...,\hat{y}(\phi(a_{m}))).
\end{align*}
for every $y\in Y$. Thus $h(a_{1},...,a_{m})=E\text{-}\sup A$.

Conversely, suppose $h(a_{1},...,a_{m})=b$ for some $b\in E$, and let $\text{\textit{\textbf{c}}}\in\Delta_{h}$. Lemma 3.7(3) implies that $\frac{\partial h(\text{\textit{\textbf{c}}})}{\partial x_{k}}\leq h(e_{k})$ for every $k\in\lbrace 1,...,m\rbrace$, and thus $\sum\limits_{k=1}^{m}\frac{\partial h(\text{\textit{\textbf{c}}})}{\partial x_{k}}a_{k}\leq\sum\limits_{k=1}^{m}h(e_{k})a_{k}$. Then $A$ is bounded above and hence $E^{\delta}$-$\sup A$ exists in $E^{\delta}$. Then $h(a_{1},..,a_{m})=E^{\delta}$-$\sup A$, and since $h(a_{1},...,a_{m})\in E$, we have $E^{\delta}\text{-}\sup A=E\text{-}\sup A$.

(2) Assume that $m\geq 2$ and that all the partial derivatives of $h$ exist and are uniformly continuous on $\lbrace s_{\theta}:\theta\in\bigcup_{n\in\mathbb{N}}P_{n}\rbrace$. It follows from \cite{Blum} that $\text{\textit{\textbf{d}}}=s_{\theta}$ for some $\theta\in[0,\pi]^{m-2}\times[0,2\pi]\ (\theta\in[0,2\pi],\ \text{for}\ m=2)$ whenever $\text{\textit{\textbf{d}}}\in\mathbb{R}^{m}$ and $||\text{\textit{\textbf{d}}}||=1$. In particular, if $\text{\textit{\textbf{d}}}\in(\mathbb{R}^{+})^{m}$ and $||\text{\textit{\textbf{d}}}||=1$ then $\text{\textit{\textbf{d}}}=s_{\theta}$ for some $\theta\in[0,\frac{\pi}{2}]^{m-1}$. Moreover, a standard induction argument verifies that $||s_{\theta}||=1$ for every $\theta\in[0,\frac{\pi}{2}]^{m-1}$. Evidently, the sequence $(\sigma_{n})$ defined by $\sigma_{n}=\sup\lbrace\nabla h(s_{\theta})\cdot\text{\textit{\textbf{a}}}:\theta\in P_{n}\rbrace$ for $n\in\mathbb{N}$ is increasing and $\sup\lbrace\sigma_{n}:n\in\mathbb{N}\rbrace=h(a_{1},...,a_{m})$. Let $r,n\in\mathbb{N}$ and assume $r<n$. For every $j\in\lbrace 1,...,m-1\rbrace$, set $I_{j}=\Bigl[\frac{(l_{j}-1)\pi}{2^{r+1}},\frac{l_{j}\pi}{2^{r+1}}\Bigr]$. By Exercise 91.10 in \cite{Zan2},
\begin{align*}
|\sup\lbrace\nabla h(s_{\theta})\cdot\text{\textit{\textbf{a}}}:\theta\in P_{n}\rbrace&-\sup\lbrace\nabla h(s_{\phi})\cdot\text{\textit{\textbf{a}}}:\phi\in P_{r}\rbrace|\\
&\leq\sup\lbrace|\nabla h(s_{\theta}-s(\phi))\cdot\text{\textit{\textbf{a}}}|:\phi\in P_{r},\theta\in\times_{j=1}^{m-1}I_{j}\rbrace\\
&\leq\sup\lbrace|\nabla h(s_{\theta}-s_{\phi})|\cdot|\text{\textit{\textbf{a}}}|:\phi\in P_{r},\theta\in\times_{j=1}^{m-1}I_{j}\rbrace.
\end{align*}
\noindent
Note that $||\phi-\theta||\leq\sqrt{m-1}\frac{\pi}{2^{r+1}}$ for every $\phi\in P_{r}$ and every $\theta\in\times_{j=1}^{m-1}I_{j}$. Thus, given $\epsilon>0$ we have for sufficiently large $r$ that
\begin{align*}
\sup\lbrace|\nabla h(s_{\theta}-s_{\phi})|\cdot|\text{\textit{\textbf{a}}}|:\ \phi\in P_{r},\theta\in\times_{j=1}^{m-1}I_{j}\rbrace\leq\epsilon\sum\limits_{k=1}^{m}|a_{k}|.
\end{align*}
\noindent
From here it is straightforward to show that $\sigma_{n}\overset{ru}{\rightarrow}h(a_{1},...,a_{m})$. The proof of the concave case is similar.
\end{proof}

\smallskip
In particular, the last part of Theorem 3.8 holds for all the \textit{$p$th power means}, where the $p$th power mean is the Stolarsky mean $\mu_{p,2p}\ (p\in\mathbb{N})$. Indeed, all of the $p$th power means are continuously differentiable on the compact set $\lbrace s_{\theta}:\theta\in[0,\frac{\pi}{2}]^{m-1}\rbrace$.

\smallskip
Two special cases of Theorem 3.8 follow as corollaries.

\smallskip
\th{Corollary 3.9.} Let $m\in\mathbb{N}\setminus\lbrace 1\rbrace$. Define $h(x_{1},...,x_{m})=\Bigl(\sum\limits_{k=1}^{m}x_{k}^{2}\Bigr)^{\frac{1}{2}}\ (x_{1},...,x_{m}\in\mathbb{R})$, and let $E$ be an $h$-complete Archimedean vector lattice. Then for every $a_{1},...,a_{m}\in E^{+}$ we have
\begin{align*}
h(a_{1},...,a_{m})=\sup\lbrace s_{\theta}(a_{1},...,a_{m}):\theta\in[0,\frac{\pi}{2}]^{m-1}\rbrace.
\end{align*}
\endth

\begin{proof} Evidently, $h$ is a member of $\mathcal{H}(\mathbb{R}^{m})$ and $h$ is convex on $(\mathbb{R}^{+})^{m}$. Also, $\Delta_{h}=\lbrace\text{\textit{\textbf{c}}}\in(\mathbb{R}^{+})^{m}: ||\text{\textit{\textbf{c}}}||=1\rbrace$. Thus for $a_{1},...,a_{m}\in E^{+}$ we have from Theorem 3.8 that
\begin{align*}
h(a_{1},...,a_{m})&=\sup\Bigl\lbrace\sum\limits_{k=1}^{m}\frac{c_{k}}{\sqrt{c_{1}^{2}+\dots+c_{m}^{2}}}a_{k}:c\in(\mathbb{R}^{+})^{m}, ||\text{\textit{\textbf{c}}}||=1\Bigr\rbrace\\
&=\sup\Bigl\lbrace\sum\limits_{k=1}^{m}d_{k}a_{k}:\text{\textit{\textbf{d}}}\in(\mathbb{R}^{+})^{m}, ||\text{\textit{\textbf{d}}}||=1\Bigr\rbrace\\
&=\sup\lbrace s_{\theta}(a_{1},...,a_{m}):\theta\in[0,\frac{\pi}{2}]^{m-1}\rbrace.
\end{align*}
\end{proof}

\smallskip
\noindent
\th{Corollary 3.10.} Define $h(x_{1},...,x_{m})=m\Bigl(\prod\limits_{k=1}^{m}|x_{k}|\Bigr)^{\frac{1}{m}}\ (x_{1},...,x_{m}\in\mathbb{R})$, and let $E$ be an $h$-complete Archimedean vector lattice. Also let $a_{1},...,a_{m}\in E^{+}$. Then
\begin{align*}
h(a_{1},...,a_{m})=\inf\lbrace\theta_{1}a_{1}+\dots+\theta_{m}a_{m}:\theta_{1},...,\theta_{m}\in(0,\infty), \theta_{1}\cdots\theta_{m}=1\rbrace.
\end{align*}
\endth

\begin{proof} Note that $h$ is a member of $\mathcal{H}(\mathbb{R}^{m})$ that is concave on $(\mathbb{R}^{+})^{m}$. Observe that $\Delta_{h}=\lbrace \text{\textit{\textbf{c}}}\in(0,\infty)^{m}: ||\text{\textit{\textbf{c}}}||=1\rbrace$. It follows from Theorem 3.8 that
\begin{align*}
h(a_{1},...,a_{m})&=\inf\Bigl\lbrace\sum\limits_{k=1}^{m}\frac{c_{1}\cdots\hat{c}_{k}\cdots c_{m}}{(c_{1}\cdots c_{m})^\frac{m-1}{m}}a_{k}:\text{\textit{\textbf{c}}}\in(0,\infty)^{m}, ||\text{\textit{\textbf{c}}}||=1\Bigr\rbrace\\
&=\inf\lbrace\theta_{1}a_{1}+\dots+\theta_{m}a_{m}:\theta_{1},...,\theta_{m}\in(0,\infty), \theta_{1}\cdots\theta_{m}=1\rbrace.
\end{align*}
\end{proof}

\smallskip
In Remark 4 of \cite{Az}, Azouzi constructs a vector lattice that he calls the \textit{square mean closure} of a given Archimedean vector lattice inside its Dedekind completion. Although Azouzi does not mention functional calculus, it turns out that his square mean closure is with respect to $\mu_{2,4}$ in Example 3.4. Indeed, Azouzi calls an Archimedean vector lattice $E$ \textit{square mean closed} if
\begin{align*}
f\boxplus g:=\sup\lbrace f\cos\theta+g\sin\theta:\theta\in[0,2\pi]\rbrace
\end{align*}
\noindent
exists in $E$ for every $f,g\in E^{+}$. Thus, for an Archimedean vector lattice $E$ we have $\mu_{2,4}(f,g)=\frac{1}{\sqrt{2}}(f\boxplus g)$ for every $f,g\in E^{+}$ (Corollary 3.9, with $m=2$). Evidently, $\mu_{2,4}(x,y)=\mu_{2,4}(|x|,|y|)\ (x,y\in\mathbb{R})$ and $f\boxplus g=|f|\boxplus|g|\ (f,g\in E)$, provided that $|f|\boxplus|g|$ exists in $E$. Thus an Archimedean vector lattice $E$ is $\mu_{2,4}$-complete if and only if $E$ is square mean closed, in which case $\mu_{2,4}(f,g)=\frac{1}{\sqrt{2}}(f\boxplus g)$ for every $f,g\in E$. Therefore, in a way, Theorem 3.8(2) generalizes the Beukers-Huijsmans-de Pagter circle approximation theorem (see section 2 of \cite{BeuHuidP}) for the existence of a modulus in the vector space complexification of a uniformly complete Archimedean vector lattice, which was later reformulated in Lemma 2.8 of \cite{Az} for square mean closed Archimedean real vector lattices. Indeed, Lemma 2.8 of \cite{Az} could be stated in terms of $\mu_{2,4}$, which is then generalized to all $\mu_{p,2p}$ in Theorem 3.8(2).

The authors of \cite{AzBoBus} call an Archimedean vector lattice $E$ \textit{geometric mean closed} if
\begin{align*}
f\boxtimes g:=\frac{1}{2}\inf\lbrace\theta f+\theta^{-1}g:\theta\in(0,\infty)\rbrace
\end{align*}
\noindent
exists in $E$ for every $f,g\in E^{+}$. For an Archimedean vector lattice $E$ we have $\mu_{1,-1}(f,g)=f\boxtimes g$ for every $f,g\in E^{+}$ (Corollary 3.10, with $m=2$).

We expand on Azouzi's idea of a square mean closure by completing Archimedean vector latices with respect to any nonempty subset of $\bigcup_{m\in\mathbb{N}}\mathcal{H}(\mathbb{R}^{m})$.

\smallskip
\noindent
\th{Definition 3.11.} For $\mathcal{D}\subseteq \bigcup_{m\in\mathbb{N}}\mathcal{H}(\mathbb{R}^{m})\ (\mathcal{D}\neq\varnothing)$ and an Archimedean vector lattice $E$, we call a pair $(E^{\mathcal{D}},\phi)$ a $\mathcal{D}$-completion of $E$ if the following hold.

\begin{itemize}
\item[(1)] $E^{\mathcal{D}}$ is a $\mathcal{D}$-complete Archimedean vector lattice.

\item[(2)] $\phi:E\rightarrow E^{\mathcal{D}}$ is an
injective vector lattice homomorphism.

\item[(3)] For every $\mathcal{D}$-complete
Archimedean vector lattice $F$ and for
every vector lattice homomorphism $T:E\rightarrow F$ there exists a
unique vector lattice homomorphism $T^{\mathcal{D}}:E^{\mathcal{D}%
}\rightarrow F$ such that $T^{\mathcal{D}}\circ\phi=T$.
\end{itemize}

\noindent
Given $h\in\mathcal{H}(\mathbb{R}^{m})$, we denote a space that satisfies (1)-(3) above for $\mathcal{D}=\lbrace h\rbrace$ by $E^{h}$, and we call $E^{h}$ an $h$-completion of $E$. We also refer to $\mathcal{D}$-completions as functional completions when the specificity of the set $\mathcal{D}$ is not present.
\endth

\smallskip
We will prove the existence and uniqueness of the $\mathcal{D}$-completion of an Archimedean vector lattice for which
we need several prerequisite results. The first of these, in a way, captures
the idea of functional calculus (see \cite{BusdPvR}) via a property of
vector lattice homomorphisms. We note that a proof of the first part of the theorem can be found in Proposition 3.6 of \cite{Kus} for uniformly complete vector lattices, and for $\mathcal{D}=\bigcup_{m\in \mathbb{N}}\mathcal{H}(\mathbb{R}^{m})$. For convenience, we write $\delta(h)$ instead of $m$ when $h\in\mathcal{H}(\mathbb{R}^{m})$.

\smallskip
\noindent
\th{Theorem 3.12.} If $E$ and $F$ are $\mathcal{D}$-complete Archimedean vector lattices, where $\mathcal{D}\subseteq\bigcup_{m\in \mathbb{N}}\mathcal{H}(\mathbb{R}^{m})$ and $\mathcal{D}$ is nonempty, and $T:E\rightarrow F$ is a vector lattice homomorphism then
\begin{equation*}
T(h(a_{1},...,a_{\delta (h)}))=h(T(a_{1}),...,T(a_{\delta (h)}))
\end{equation*}
\noindent
for every $h\in \mathcal{D}$ and for every $a_{1},...,a_{\delta (h)}\in E$. Moreover, if there exists $h\in 
\mathcal{D}$ such that $h(\epsilon _{1}x,...,\epsilon _{\delta
(h)}x)=\lambda |x|\ (x\in \mathbb{R})$ for some $\epsilon
_{1},...,\epsilon _{\delta (h)}\in\mathbb{R}$ and some $\lambda\in 
\mathbb{R}\setminus \{0\}$ then every linear map $S:E\rightarrow F$ such that $S(h(a_{1},...,a_{\delta(h)}))=h(S(a_{1}),...,S(a_{\delta (h)}))$ for every $a_{1},...,a_{\delta (h)}\in E$ is a vector lattice homomorphism.
\endth

\begin{proof}
Let $\mathcal{D}\subseteq\bigcup_{m\in \mathbb{N}}\mathcal{H}(\mathbb{R}^{m})$, let $E$ and $F$ be $\mathcal{D}$-complete Archimedean vector lattices, and let $T:E\rightarrow F$ be a vector lattice homomorphism. Let $h\in \mathcal{D}$ and let $a_{1},...,a_{\delta (h)}\in E$. Define $G_{1}:=\left\langle
a_{1},...,a_{\delta (h)},h(a_{1},...,a_{\delta (h)})\right\rangle $ and $G_{2}:=\left\langle T(a_{1}),...,T(a_{\delta (h)}),T(h(a_{1},...,a_{\delta
(h)}))\right\rangle $. Let $\omega\in H(G_{2})$. Since $E$ and $F$ are $h$-complete and $\omega\circ T|_{G_{1}}\in H(G_{1})$, we have 
\begin{equation*}
h(\omega (T(a_{1})),...,\omega (T(a_{\delta (h)})))=\omega\circ
T(h(a_{1},...,a_{\delta (h)})).
\end{equation*}
\noindent
Then 
\begin{equation*}
h(T(a_{1}),...,T(a_{\delta (h)}))=T(h(a_{1},...,a_{\delta (h)})).
\end{equation*}
\noindent
Conversely, suppose that there exist $\epsilon_{1},...,\epsilon_{\delta(h)}\in\mathbb{R}$ and $\lambda\in\mathbb{R}\setminus\{0\}$ such that for every $x\in\mathbb{R}$ we have that $h(\epsilon_{1}x,...,\epsilon_{\delta(h)}x)=\lambda |x|$. Let $S:E\rightarrow F$ be a linear map such that $S(h(a_{1},...,a_{\delta(h)}))=h(S(a_{1}),...,S(a_{\delta (h)}))$ for every $a_{1},...,a_{\delta(h)}\in E$, and let $a\in E$. Also let $\omega\in H(\langle a\rangle)$. Then $h(\omega(\epsilon_{1}a),...,\omega(\epsilon_{m}a))=h(\epsilon_{1}\omega(a),...,\epsilon_{m}\omega(a))=\lambda|\omega(a)|=\omega(\lambda|a|)$, and thus $h(\epsilon _{1}a,...,\epsilon _{\delta (h)}a)=\lambda|a|$. Similarly, $h(S(\epsilon_{1}a),...,S(\epsilon _{\delta (h)}a))=\lambda |S(a)|$. Hence
\begin{equation*}
S(\lambda |a|)=S(h(\epsilon_{1}a,...,\epsilon_{\delta (h)}a))=h(S(\epsilon_{1}a),...,S(\epsilon_{\delta (h)}a))=\lambda|S(a)|\text{.}
\end{equation*}
\noindent
Then $S(|a|)=|S(a)|$ since $\lambda\neq 0$ and $S$ is linear.
\end{proof}

\smallskip
As a particular case of theorem above, suppose that $E$ and $F$ are $h$-complete Archimedean vector lattices for some Stolarsky mean or Gini mean $h$. Then a linear map $T:E\rightarrow F$ is a vector lattice homomorphism if and only if $T(h(f,g))=h(T(f),T(g))$ for every $f,g\in E$. Thus Theorem 3.12 generalizes a result by Azouzi, Boulabiar, and the first author in \cite{AzBoBus}, as well as a proposition of Azouzi in \cite{Az}. We point out that Corollary 3.13 below is a generalization of Lemma 4.3 in \cite{AzBoBus} and corrects a mistake (first noted in \cite{Chen}) in its
proof. For the proof of Corollary 3.13, respectively Corollary 3.14, apply Corollary 3.9, respectively Corollary 3.10, and Theorem 3.12.

\smallskip
\noindent
\th{Corollary 3.13.}\ \textnormal{(\cite{AzBoBus}, Corollary 4.7)} For square mean complete Archimedean vector
lattices $E$ and $F$ and a linear map $T:E\rightarrow F$, the following are equivalent.
\begin{itemize}
\item[(1)] $T$ is a vector lattice homomorphism.
\item[(2)] $T(f\boxplus g)=T(f)\boxplus T(g)$ for every $f,g\in E^{+}$.
\end{itemize}
\endth

\smallskip
\noindent
\th{Corollary 3.14.}\ \textnormal{(\cite{Az}, Proposition 2.20)} For geometric mean complete Archimedean vector lattices $E$ and $F$ and a linear map $T:E\rightarrow F$, the following are equivalent.
\begin{itemize}
\item[(1)] $T$ is a vector lattice homomorphism.
\item[(2)] $T(f\boxtimes g)=T(f)\boxtimes T(g)$ for every $f,g\in E^{+}$.
\end{itemize}
\endth

\smallskip
The following theorem is needed for our construction of functional completions.

\smallskip
\noindent
\th{Theorem 3.15. \textnormal{(\cite{BusdPvR}, Theorem 3.7)}} A uniformly complete Archimedean vector lattice $E$ is $\bigcup_{m\in\mathbb{N}}\mathcal{H}(\mathbb{R}^{m})$-complete.
\endth

\smallskip
We remind the reader that $h\in\mathcal{H}(\mathbb{R}^{m})$ is \textit{positive} if $h(x_{1},...,x_{m})\in 
\mathbb{R}^{+}$ for every $x_{1},...,x_{m}\in\mathbb{R}^{+}$. If $h(x_{1},...,x_{m})=h(|x_{1}|,...,|x_{m}|)$ for every $x_{1},...,x_{m}\in 
\mathbb{R}$, we call $h$ \textit{absolutely invariant}. We denote the set of
all $h\in\mathcal{H}(\mathbb{R}^{m})$ that are positive and absolutely invariant by $\mathcal{H}_{|\ |}^{+}(\mathbb{R}^{m})$. Examples of such functions include the Stolarsky and Gini means from Examples 3.4 and 3.5. We first manufacture a $\mathcal{D}$-completion $E^{\mathcal{D}}$ of $E$ for any Archimedean vector lattice $E$. That this is indeed the $\mathcal{D}$-completion will subsequently be proved. Thus let $E$ be an Archimedean vector lattice and assume that $\mathcal{D}\subseteq\bigcup_{m\in\mathbb{N}}\mathcal{H}(\mathbb{R}^{m})$ is nonempty. Let $(E^{u},\phi)$ be the uniform completion of $E$. Following the lead of Azouzi in Remark 4 of \cite{Az}, define
\begin{equation*}
\begin{array}{l}
E_{1}:=\phi (E),\ \text{and for every}\ n\in\mathbb{N}, \\
E_{n+1}:=\left\langle E_{n}\cup \{h(a_{1},...,a_{\delta (h)}):h\in \mathcal{D},a_{1},...,a_{\delta (h)}\in E_{n}\}\right\rangle,
\end{array}
\end{equation*}
\noindent
where the latter is the vector lattice generated in $E^{u}$. We define 

\begin{equation*}
E^{\mathcal{D}}:=\bigcup_{n\in \mathbb{N}}E_{n}.
\end{equation*}

\noindent
The following proposition is immediate.

\smallskip
\noindent
\th{Proposition 3.16.} $E^{\mathcal{D}}$ is $\mathcal{D}$-complete.
\endth

\smallskip
By using Proposition 3.2(2) one can, alternatively to the definition of $\mathcal{D}$-completion, replace the homomorphisms in that definition by positive maps if the range space is required to be uniformly complete. This is the content of the next proposition.

\smallskip
\noindent
\th{Proposition 3.17.} Let $E_{1},...,E_{s},F$ be Archimedean vector lattices with $F$ uniformly complete. Let $\mathcal{D}$ be a nonempty subset of $\bigcup_{m\in \mathbb{N}}\mathcal{H}(\mathbb{R}^{m})$. If $T:\times_{k=1}^{s}E_{k}\rightarrow F$ is a positive $s$-linear map then there exists a unique positive $s$-linear map $T^{\mathcal{D}}:\times_{k=1}^{s}E_{k}^{\mathcal{D}}\rightarrow F$ such that $T^{\mathcal{D}}\circ\phi=T$.
\endth

\smallskip
We prove that $E^{\mathcal{D}}$ is the $\mathcal{D}$-completion of $E$ by proving more, involving multilinear maps, as follows.

\smallskip
\noindent
\th{Theorem 3.18.} Let $E_{1},...,E_{s},F$ be Archimedean vector lattices and assume that $F$ is $\mathcal{D}$-complete \textnormal{(}$\mathcal{D}\subseteq\bigcup_{m\in\mathbb{N}}\mathcal{H}(\mathbb{R}^{m}),\ \mathcal{D}\neq\varnothing$\textnormal{)}. Also let $T:\times_{k=1}^{s}E_{k}\rightarrow F$ be a vector lattice $s$-morphism. Denoting, for every $k\in \{1,...,s\}$, the natural embedding of $E_{k}$ into $E_{k}^{u}$ by $\phi_{k}$, the following hold.
\begin{itemize}
\item[(1)] If all elements of $\mathcal{D}$ are absolutely invariant then there exists a unique vector lattice $s$-homomorphism $T^{\mathcal{D}}:\times _{k=1}^{s}E_{k}^{\mathcal{D}}\rightarrow F$ such that 
\begin{equation*}
T^{\mathcal{D}}(\phi_{1}(f_{1}),...,\phi_{s}(f_{s}))=T(f_{1},...,f_{s})
\end{equation*}
\noindent
for every $f_{k}\in E_{k}\ (k\in \{1,...,s\})$.
\item[(2)] If $s=1$ then statement \textnormal{(1)} holds for all nonempty $\mathcal{D}\subseteq\bigcup_{m\in\mathbb{N}}\mathcal{H}(\mathbb{R}^{m})$.
\end{itemize}
\endth

\begin{proof}
We will prove statement (1) and (2) simultaneously, since the absolute invariance of elements of $\mathcal{D}$ is not
used in our proof of (1) in case $s=1$. Let $\mathcal{D}\subseteq\bigcup_{m\in\mathbb{N}}\mathcal{H}_{|\ |}^{+}(\mathbb{R}^{m})$ be nonempty and let $E_{1},...,E_{s},F$ be Archimedean vector lattices. Let $T:\times_{k=1}^{s}E_{k}\rightarrow F$ be a vector lattice $s$-morphism and assume $j\in \{1,...,s\}$. Then $E_{j}^{\mathcal{D}}$, as defined preceding Proposition 3.16, is a $\mathcal{D}$-complete Archimedean vector lattice and we denote the
natural embedding into the uniform completion of $E_{j}$ by $\phi_{j}$. Clearly $\phi_{j}(E_{j})\subseteq E_{j}^{\mathcal{D}}$ and $\phi_{j}$ is an injective vector lattice homomorphism from $E_{j}$ into $E_{j}^{\mathcal{D}}$. We
prove that there exists a unique vector lattice $s$-morphism $T^{\mathcal{D}}:\times_{k=1}^{s}E_{k}^{\mathcal{D}}\rightarrow F$ such that $T^{\mathcal{D}}(\phi_{1}(f_{1}),...,\phi_{s}(f_{s}))=T(f_{1},...,f_{s})$ for every $f_{k}\in E_{k}\ (k\in \{1,...,s\})$. We consider $T$ to be a vector lattice $s$-morphism from $\times_{k=1}^{s}E_{k}$ to $F^{u}$. By Proposition 3.2(2) there exists a unique vector lattice $s$-morphism $T^{u}:\times_{k=1}^{s}E_{k}^{u}\rightarrow F^{u}$ such that $T^{u}(\phi_{1}(f_{1}),...,\phi_{s}(f_{s}))=T(f_{1},...,f_{s})$ for every $f_{k}\in E_{k}\ (k\in \{1,...,s\})$. Define $T^{\mathcal{D}}:=T^{u}|_{\times_{k=1}^{s}E_{k}^{\mathcal{D}}}$. To prove that $T^{\mathcal{D}}(\times_{k=1}^{s}E_{k}^{\mathcal{D}})\subseteq F$, we write $E_{k,n}$ for $(E_{k})_{n}$, where $(E_{k})_{n}$ is defined as preceding Proposition 3.16, and we use induction with respect to $n$. Again for convenience, we write $\delta(h)$ instead of $m$ when $h\in\mathcal{H}(\mathbb{R}^{m})$.

Obviously $T^{\mathcal{D}}(\times_{k=1}^{s}E_{k,1})\subseteq F$. Let $n\in\mathbb{N}$ and suppose that $T^{\mathcal{D}}(\times_{k=1}^{s}E_{k,n})\subseteq F$. Let $h_{1},...,h_{s}\in \mathcal{D}$ and let $a_{1}^{k},...,a_{\delta (h_{k})}^{k}\in E_{k,n}\ (k\in \{1,...,s\})$. Write
\begin{equation*}
x=T^{\mathcal{D}}\bigl(h_{1}(a_{1}^{1},...,a_{\delta
(h_{1})}^{1}),...,h_{s}(a_{1}^{s},...,a_{\delta (h_{s})}^{s})\bigr).
\end{equation*}
\noindent
Since $h_{k}$ is absolutely invariant for each $k\in\lbrace 1,...,s\rbrace$ we may assume that $a_{1}^{k},...,a_{\delta (h_{k})}^{k}\in E_{k,n}^{+}\ (k\in
\{1,...,s\}) $. Given an Archimedean vector lattice $E$ and a nonempty subset $A$ of $E$, and $h\in\mathcal{H}(\mathbb{R}^{m})$, define $h(A):=\{f\in E:f=h(a_{1},...,a_{m})\ \text{for some}\ a_{1},...,a_{m}\in A\}$. Then $x\in T^{\mathcal{D}}\bigl(h_{1}(E_{1,n}^{+})\times\dots\times h_{s}(E_{s,n}^{+})\bigr)$. By repeatedly employing Theorem 3.12, we have
\begin{gather*}
x\in h_{1}\Bigl(\bigl(T^{\mathcal{D}}(E_{1,n}^{+}\times
h_{2}(E_{2,n}^{+})\times\dots\times h_{s}(E_{s,n}^{+}))\bigr)^{\delta
(h_{1})}\Bigr)\\
\subseteq h_{1}\Biggl(\biggl(h_{2}\Bigl(\bigl(T^{\mathcal{D}}(E_{1,n}^{+}\times
E_{2,n}^{+}\times h_{3}(E_{3,n}^{+})\times\dots\times h_{s}(E_{s,n}^{+}))\bigr)^{\delta (h_{2})}\Bigr)\biggr)^{\delta (h_{1})}\Biggr) \\
\subseteq h_{1}\Biggl(\Biggl(\dots h_{s-1}\Biggl(\biggl(h_{s}\Bigl(\bigl(T^{%
\mathcal{D}}(E_{1,n}^{+}\times \dots \times E_{s,n}^{+})\bigr)^{\delta
(h_{s})}\Bigr)\biggr)^{\delta (h_{s-1})}\Biggr)\dots \Biggr)^{\delta (h_{1})}%
\Biggr) \\
\subseteq F,
\end{gather*}
\noindent
where the last inclusion follows from the induction
hypothesis and the assumption that $F$ is $\mathcal{D}$-complete. Moreover,
from the $s$-linearity of $T^{\mathcal{D}}$,
\begin{equation*}
T^{\mathcal{D}}(\times _{k=1}^{s}[E_{k,n}\bigcup \{h(a_{1}^{k},...,a_{\delta
(h)}^{k}):h\in \mathcal{D},a_{1}^{k},...,a_{\delta (h)}^{k}\in
E_{k,n}\}])\subseteq F.
\end{equation*}
\noindent
By Exercise 4.1.8 in \cite{AB}, every element of $E_{k,n+1}^{+}$ can be expressed as $\bigwedge\limits_{j=1}^{p_{k}}\bigvee\limits_{l=1}^{q_{k}}u_{k,j,l}$ for some $u_{k,j,l}\in[E_{k,n}\bigcup\{h(a_{1}^{k},...,a_{\delta (h)}^{k}):h\in \mathcal{D},a_{1}^{k},...,a_{\delta (h)}^{k}\in E_{k,n}\}]$. Moreover, we may assume that each $u_{k,j,l}$ is positive. Since $T^{\mathcal{D}}$
is a vector lattice $s$-morphism,
\begin{equation*}
T^{\mathcal{D}}(\bigwedge\limits_{j=1}^{p_{1}}\bigvee%
\limits_{l=1}^{q_{1}}u_{1,j,l},...,\bigwedge\limits_{j=1}^{p_{s}}\bigvee%
\limits_{l=1}^{q_{s}}u_{s,j,l})=\bigwedge\limits_{j=1}^{p_{1}}\bigvee%
\limits_{l=1}^{q_{1}}\dots
\bigwedge\limits_{j=1}^{p_{s}}\bigvee\limits_{l=1}^{q_{s}}T^{\mathcal{D}%
}(u_{1,j,l},...,u_{s,j,l})
\end{equation*}
\noindent
for every $u_{k,j,l}\in E_{k,n+1}^{+}$. Since $T^{\mathcal{D}}(u_{1,j,l},...,u_{s,j,l})\in F$ for each $j$ and each $l$, it
follows that
\begin{equation*}
\bigwedge\limits_{j=1}^{p_{1}}\bigvee\limits_{l=1}^{q_{1}}\dots
\bigwedge\limits_{j=1}^{p_{s}}\bigvee\limits_{l=1}^{q_{s}}T^{\mathcal{D}%
}(u_{1,j,l},...,u_{s,j,l})\in F,
\end{equation*}
\noindent and hence $T^{\mathcal{D}}(\times
_{k=1}^{s}E_{k,n+1}^{+})\subseteq F$. Then also $T^{\mathcal{D}}(\times
_{k=1}^{s}E_{k,n+1})\subseteq F$ because $T^{\mathcal{D}}$ is $s$-linear. This completes the proof.
\end{proof}

\smallskip
\noindent
\th{Corollary 3.19.} Let $E$ be an Archimedean vector lattice and let $\mathcal{D}\subseteq\bigcup_{m\in\mathbb{N}}\mathcal{H}(\mathbb{R}^{m})$ be nonempty. Then $E^{\mathcal{D}}$ with the natural embedding from $E$ into $E^{\mathcal{D}}$ is the unique $\mathcal{D}$-completion of $E$.
\endth

\begin{proof}
Let $E$ be an Archimedean
vector lattice and let $\mathcal{D}\subseteq\bigcup_{m\in\mathbb{N}}\mathcal{H}(\mathbb{R}^{m})\ (\mathcal{D}\neq\varnothing)$. We proved in the previous theorem that $E^{\mathcal{D}}$ with the natural embedding from $E$ into $E^{\mathcal{D}}$ is a $\mathcal{D}$-completion of $E$. Next, we prove the uniqueness of $E^{\mathcal{D}}$. Suppose $(E_{1}^{\mathcal{D}},\phi _{1})$ and $(E_{2}^{\mathcal{D}},\phi_{2})$ are $\mathcal{D}$-completions of $E$. Since $\phi_{1}:E\rightarrow E_{1}^{\mathcal{D}}$ is a vector lattice homomorphism, there exists a unique
vector lattice homomorphism $\phi_{1}^{\mathcal{D}}:E_{2}^{\mathcal{D}}\rightarrow E_{1}^{\mathcal{D}}$ such that $\phi_{1}^{\mathcal{D}}\circ\phi_{2}=\phi_{1}$. Likewise, there exists a unique vector lattice
homomorphism $\phi_{2}^{\mathcal{D}}:E_{1}^{\mathcal{D}}\rightarrow E_{2}^{\mathcal{D}}$ such that $\phi_{2}^{\mathcal{D}}\circ\phi_{1}=\phi_{2}$.
Then we have $\phi_{2}^{\mathcal{D}}\circ\phi_{1}^{\mathcal{D}}\circ\phi_{2}=\phi_{2}^{\mathcal{D}}\circ\phi_{1}=\phi_{2}$. Thus $\phi_{2}^{\mathcal{D}}\circ\phi_{1}^{\mathcal{D}}=I$. Similarly, $\phi_{1}^{\mathcal{D}}\circ\phi_{2}^{\mathcal{D}}=I$.
\end{proof}

\end{document}